\documentclass[a4paper, 12pt]{amsart}
\usepackage{amssymb,amsmath,bm}
\usepackage{color}
\usepackage{amsthm}
\definecolor{light-red}{rgb}{0.85,0,0}%

\long\def\ori#1{\relax\ignorespaces}

\def\dotbigcup{\mathbin{\dot{\bigcup}}}

\def\mathbf#1{{\bm #1}}
\newtheorem{thm}{Theorem}[section]
\newtheorem{cor}[thm]{Corollary}
\newtheorem{lemma}[thm]{Lemma}

\theoremstyle{definition}
\newtheorem{defn}[thm]{Definition}
\theoremstyle{remark}
\newtheorem{rem}[thm]{Remark}

\numberwithin{equation}{section}
\newcommand{\GS}{Gr\"obner-Shirshov }
\newcommand{\vb}{\mathord{\hskip 1pt|\hskip 1pt}}

\newcommand{\CD}{Composition-Diamond}
\begin{document}

\title[Markov and Artin normal form theorem for braid groups]
{Markov and Artin normal form theorem for braid groups}

\author[Bokut]{L. A. Bokut$^*$}
\address{Sobolev Institute of Mathematics, Novosibirsk 630090, Russia}
\email{bokut@math.nsc.ru}
\thanks{$^*$Supported in part by the Russia's Fund for
    Fundamental Research and the Leading Scientific Schools Fund (Russia).}

\author[Chaynikov]{V.V. Chaynikov$^*$}
\address{Faculty of mathematics and mechanics, Novosibirsk State University, 
Novosibirsk, 630090, Russia}
\email{vvch@math.nsc.ru}

\author[Shum]{K.P. Shum$^{**} $ }
\address{Faculty of Science, The Chiness University of Hong Kong, Hong Kong, China (SAR)}
\email{kpshum@math.cuhk.edu.hk}
\thanks{$^{**}$Partially supported by a RGC(HK) grant 2003/04}

\maketitle
\begin{abstract}
In this paper we will present the results of Artin--Markov on braid groups by using the \GS basis.
As a consequence we can reobtain the normal form of Artin--Markov--Ivanovsky as an easy corollary.
\end{abstract}
\section{Introduction}
The theory of braid groups were first studied by E.
Artin dated back to 1925 (see ~\cite{Ar25}). Artin
established generators and defining relations of the braid group and
 a faithful representation theorem of the braid group as
a subgroup of the automorphism group  of a free group. Later on, 
Markov ~\cite{Mar45} (1945) and Artin ~\cite{Ar47} (1947) enriched the Artin generators
by the Burau elements (~\cite{Bur32}) and find a presentation of the braid group in the
Artin-Burau generators. We call it the Artin-Markov presentation of the braid group in 
the Artin-Burau generators.
The main result of (~\cite{Mar45}) and (~\cite{Ar47})
is the normal form theorem for the braid group. We call it the Artin-Markov-Ivanovsky normal form 
theorem for Markov credited the result to A. Ivanovsky; it looks that there is no paper 
published by A. Ivanovskii himself.
 However, the proof of the
normal form  theorem of the braid group given by Markov and Artin are
rather complicated because the automorphism of the free groups are
involved. In this note we will try to escape the technical details
of the proof given by Markov ~\cite{Mar45} and Artin ~\cite{Ar25}.
We will show by  direct calculations (of compositions) that 
the Artin-Markov presentation of the braid group is the
minimal \GS basis of the braid group in the Artin-Burau
generators and under an appropriate ordering of group words,
namely, the inverse tower ordering. As a consequence, by the
Composition-Diamond Lemma the Artin-Markov-Ivanovsky normal form in the
braid group follows as an immediate corollary of our result.

\section{Basic notations and results}
In this section we give some basic notations and cite some useful
results in the literature. We first let  $\mathbf{X}$ be a
linearly ordered set and  $k$ a field. Let $k\langle
\mathbf{X}\rangle$ the free associative algebra over $\mathbf{X}$
and $k$. On the set $\mathbf{X}^*$ of words generated by
$\mathbf{X}$ we impose a well order $\leq$ compatible with the
multiplication of words. We call this kind of order a monomial order.

Now, let $f\in k\langle \mathbf{X}\rangle$ be a polynomial with
leading word $\bar f$.  We say that $f$ is monic if $\bar f$
occurs in $f$ with coefficient $1$. We now formulate the following
definitions.
\begin{defn}
Let $f$ and $g$ be two monic polynomials.
\begin{itemize}
\item[(i)]
   Let $w$ be a word such that $w = \bar fb= a \bar g$, with $\deg(
   \bar f) + \deg(\bar g) > \deg (w)$. Then we call the polynomial
   $(f,g)_w$ the intersection composition of $f$ and $g$ with respect
   to $w$ if $(f,g)_w = fb - ag$.
\item[(ii)]
   Let $w = \bar f = a\bar gb$. Then we call the polynomial
   $(f,g)_w = f - agb$ the inclusion composition of $f$ and $g$ with respect
   to $w$.

   In the above case, the transformation
   $ f\mapsto (f,g)_w = f -agb$
   is called the elimination of the leading word (ELW) of g in f.

\item[(iii)]
   Let $\mathbf{S} \subset k\langle \mathbf{X}\rangle$. We call a 
composition $(f, g)_w$ trivial relative to $\mathbf{S}$ (and $w$) if
   $$
            (f,g)_w = \sum \alpha_i a_it_i b_i,
   $$
   where   $t_i \in \mathbf{S}$, $a_i, b_i \in  \mathbf{X}  ^*$, and
   $\bar {a_it_ib_i} < w$. In a notation,  $(f, g)_w\equiv
   0\bmod (\mathbf{S}, w)$.
   In particular, if $(f,g)_w$ goes to zero
   by using ELW of polynomials from $\mathbf{S}$,
   then $(f,g)_w$ is trivial relative to $\mathbf{S}$. 
 We assume that  $f_1$ and $f_2$ are
some polynomials satisfying the condition $ f_1 \equiv f_2 \bmod
(\mathbf{S}, w) $ if $ f_1-f_2 \equiv 0 \bmod (\mathbf{S},w)$.
\end{itemize}
\end{defn}

\begin{defn}
Then we call $\mathbf{S}$ a \emph{\GS\  set (basis)} in $k\langle \mathbf{X}\rangle$ 
if any composition
of polynomials from $\mathbf{S}$ is trivial relative to $\mathbf{S}$.
\end{defn}

 From now on,
we denote the algebra with set of generators $\mathbf{X}$ and set of defining 
relations  $\mathbf{S}$
by  $\langle \mathbf{X}\vb \mathbf{S}\rangle$.

 The following lemma and its applications to \GS\ bases was due to
 Newman ~\cite{Ne42}, Shirshov
(~\cite{Sh62}), Buchberger ~\cite{Bu65}, ~\cite{Bu70} and Bergman
~\cite{Be78}. The lemma given below was formulated by Bokut (see
 ~\cite{Bo72}, ~\cite{Bo76}).

\begin{lemma} [Composition-Diamond Lemma] Let $R=\langle \mathbf{X}\vb \mathbf{S}\rangle$.
The set of defining relations   $\mathbf{S}$ is a \GS\ set  if
and only if the set
$$
    Irr(S)=   \{u \in \mathbf{X}^* \mid u \neq a\bar fb, \mbox{ for any }f\in \mathbf{S}\}
$$
 of $\mathbf{S}$-irreducible words consists of  a linear basis of $R$.
\end{lemma}

\begin{defn} Let $\mathbf{S}$  be a \GS basis in   $k\langle
\mathbf{X}\rangle$. Then  $\mathbf{S}$ is called a  \emph{minimal}
\GS\ basis if for any $s\in \mathbf{S}$, $s$ is a linear
combination of $\mathbf{S}\setminus\{s\}$-irreducible words. Any ideal
of $k\langle \mathbf{X} \rangle$ has a unique minimal \GS basis 
(i.e., a set of generators of the ideal).
\end{defn}

If a subset $\mathbf{S}$ of $k\langle \mathbf{X}\rangle$ is not a
\GS\ basis then one can add to $\mathbf{S}$ all nontrivial
compositions of polynomials of
 $\mathbf{S}$, and continue this process (infinitely) many times  in 
order to have a
\GS\  set (basis) $\mathbf{S}^{\rm comp}$ that generates the same ideal as $S$. 
This procedure is called the Buchberger - Shirshov algorithm
\cite{Sh62}, \cite{Bu65}, \cite{Bu70}.
 
\begin{defn}
A polynomial  $f$ is called semigroup relation if $f$ is of the
form $u - v$, where $u,v \in \mathbf{X}^*$.  If 
$S$ is a set of semigroup relations,  then any nontrivial
composition of $S$  has the same form. Consequently, the set
$\mathbf{S}^{\rm comp}$ also consists of semigroup relations.
\end{defn}

\begin{rem}
Let $A = sgp\langle \mathbf{X}\vb \mathbf{S}\rangle$ be a
semigroup presentation. Then $\mathbf{S} \subset k\langle
\mathbf{S}\rangle$ and one can find a \GS\ basis $
\mathbf{S}^{\rm comp}$. This set does not depend on $k$ and it
consists of semigroup relations. In this case, we call
$\mathbf{S}^{\rm comp}$ a \GS basis of $A$. It is the same as a
\GS\ basis of the semigroup algebra $kA = \langle
\mathbf{X}\vb\mathbf{S}\rangle$.
\end{rem}

We now introduce the concept of inverse tower ordering of words.
\begin{defn}
Let $X=Y\dot{\bigcup}Z$, words $Y^*$ and the letters $Z$ are well
ordered. Suppose that the order on $Y^*$ is monomial. Any word in
$X$ has the form $ u=u_0z_1\dots u_{k-1}z_ku_{k}, $  where $k\geq
0$, $z_i\in Z$, $u_i\in Y^*$. Define the inverse weight of the word
$u \in X$ by:
$$
inwt(u)=(k,u_k,z_k,\dots ,z_1,u_0).
$$
Now we order inverse weights lexicographically and define
$$
u>v \Longleftrightarrow inwt(u)>inwt(v).
$$
Then we call the above order the \emph{inverse tower order}. Clearly,
the above  order is a monomial order.

In case $Y=T\dotbigcup U$ and $Y^*$ are endowed with
the inverse tower order, we call
order of words on $X$ the inverse tower order of word
relative the presentation
$$
X=(T\dotbigcup U)\dotbigcup Z.
$$

In general, we can define the inverse tower order of $X$-words
relative to the presentation
$$
X=(\dots (X^{(n)}\dotbigcup X^{(n-1)})\dotbigcup \dots )\dotbigcup
X^{(0)},
$$
where $X^{(n)}$-words are endowed by a monomial order.
\end{defn}

\begin{defn}
Let $ \Sigma=\{\sigma_1, \dots, \sigma_{n-1}\}$ be a finite
alphabet. Then, the following group presentation define the $n$-strand
braid group:
$$
B_n=\langle\Sigma \mid
\sigma_{i+1}\sigma_i\sigma_{i+1}=\sigma_i\sigma_{i+1}\sigma_i,
\sigma_i\sigma_j=\sigma_j\sigma_i, i-j>1\rangle
$$
 Here  any index falls into the
interval $[1,n-1]$ .

In the braid group $B_n$, we now introduce a new set of generators.
We call them the Artin-Burau generators.

In the braid group $B_n$, we let
$$
s_{i,i+1}=\sigma_i^2,\  s_{i,j+1}=\sigma_j\dots
\sigma_{i+1}\sigma_i^2\sigma_{i+1}^{-1} \dots \sigma_j^{-1},
$$
where $1\leq i<j\leq n-1$.

Form the set  $$S_j=\{s_{i,j}, s_{i,j}^{-1},
 2< i < j  <n \}$$
 and
$$\Sigma^{-1}=\{\sigma_1^{-1},\dots , \sigma_{n-1}^{-1}\}.$$
\end{defn}

Then the set
$$
S=S_n\cup S_{n-1}\cup \dots \cup S_2\cup \Sigma^{-1}
$$
generates $B_n$ as a semigroup. We call elements of $S$ the Artin-Burau generators of 
$B_n$. Observe that generators $\sigma_i$
are omitted as well as the trivial group relations on them. With
the above notation Markov ~\cite{Mar45} used
$s_{i,i+1}\sigma_i^{-1}$ to replace $\sigma_i$, and
$\sigma_i^{-2}=s_{i,i+1}^{-1}$ to replace $\sigma_i^{-1}\sigma_i=1$.

Then we order the set $S$ in the following way:

$$
S_n<S_{n-1}< \dots <S_2<\Sigma^{-1},
$$
Clearly, in the above chain, any letter of $S_n$ is less than any
letter of $S_{n-1}$ and so on. Also we define for each $j$

$$
s_{1,j}^{-1}<s_{1,j}<s_{2,j}^{-1}< \dots <s_{j-1,j}, \hbox{\quad
and\quad}\sigma_1^{-1}<\sigma_2^{-1}< \dots \sigma_{n-1}^{-1}.
$$
With above notation, we now able to order the $S$-words by using the
inverse tower order, according to the fixed presentation of $S$ as
the union of $S_j$ and $\Sigma^{-1}$. We order the $S_n$-words by the 
$deg-inlex$ order, i.e., we first compare the words by length and
than by inverse lexicographical order, starting from their last
letters.

The following abbreviations are taken from ~\cite{Mar45}.
$$
\sigma_{i,j+1}=\sigma_i^{-1}\dots \sigma_{j}^{-1}, 1\leq i\leq j\leq
n-1, \sigma_{ii}=1.
$$
Also we denote $ \{a,b\}=b^{-1}ab $.

\section{Main results}

We first cite some crucial results from ~\cite{Mar45},
~\cite{Ar47}.

The first lemma is fairly easy.  

\begin{lemma} (~\cite[Lemma ~3]{Mar45},
~\cite[p.119]{Ar47})\label{L1}
The following relations hold in the braid group $B_n$ for
($\delta=\pm1$):
\begin{align}
&\sigma_k^{-1}s_{i,j}^{\delta}=s_{i,j}^{\delta}\sigma_k^{-1},\quad
k\neq i-1,i,j-1,j,
\label{E1}\\
&\sigma_i^{-1}s_{i,i+1}^{\delta}=s_{i,i+1}^{\delta}\sigma_i^{-1},\label{E2}\\
&\sigma_{i-1}^{-1}s_{i,j}^{\delta}=s_{i-1,j}^{\delta}\sigma_{i-1}^{-1},\label{E3}\\
&\sigma_{i}^{-1}s_{i,j}^{\delta}=\{s_{i+1,j}^{\delta},s_{i,i+1}\}\sigma_{i}^{-1},\label{E4}\\
&\sigma_{j-1}^{-1}s_{i,j}^{\delta}=s_{i,j-1}^{\delta}\sigma_{j-1}^{-1},\label{E5}\\
&\sigma_{j}^{-1}s_{i,j}^{\delta}=\{s_{i,j+1}^{\delta},s_{j,j+1}\}\sigma_{j}^{-1}.\label{E6}\\
\noalign{\hbox{Also we have}}
&\sigma_{i-1}s_{i,j}^{\delta}=\{s_{i-1,j}^{\delta},s_{i-1,i}^{-1}\}\sigma_{i-1},\nonumber\\
&\sigma_{i}s_{i,j}^{\delta}=s_{i+1,j}^{\delta}\sigma_{i},\nonumber\\
&\sigma_{j-1}s_{i,j}^{\delta}=\{s_{i,j-1}^{\delta},s_{j-1,j}^{-1}\}\sigma_{j-1},\nonumber\\
&\sigma_{j}s_{i,j}^{\delta}=s_{i,j+1}^{\delta}\sigma_{j}.\nonumber
\end{align}
\end{lemma}

\begin{rem}\label{R1}
We note that last relations from the above lemma will not in \GS
basis of the braid group $B_n$ (see below). We use them for a sketch of a proof of next 
Lemma \ref{L2}. 

\end{rem}

\begin{lemma}(~\cite[Lemmas 6 and 7]{Mar45}, ~\cite[Theorem~18]{Ar47})\label{L2} 
The following relations hold in the braid group $B_n$ for all $i<j<k<l$,
$\varepsilon=\pm1$:
\begin{align}
&s_{j,k}^{-1}s_{k,l}^{\varepsilon}=\{s_{k,l}^{\varepsilon}, s_{j,l}^{-1}\}s_{j,k}^{-1},\label{E7}\\
&s_{j,k}s_{k,l}^{\varepsilon}=\{s_{k,l}^{\varepsilon}, s_{j,l}s_{k,l}\}s_{j,k},\label{E8}\\
&s_{j,k}^{-1}s_{j,l}^{\varepsilon}=\{s_{j,l}^{\varepsilon},
s_{k,l}^{-1}s_{j,l}^{-1}\}s_{j,k}^{-1},
\label{E9}\\
&s_{j,k}s_{j,l}^{\varepsilon}=\{s_{j,l}^{\varepsilon},
s_{k,l}\}s_{j,k},
\label{E10}\\
&s_{i,k}^{-1}s_{j,l}^{\varepsilon}=\{s_{j,l}^{\varepsilon},s_{k,l}s_{i,l}
s_{k,l}^{-1}s_{i,l}^{-1}\}
s_{i,k}^{-1},\label{E11}\\
&s_{i,k}s_{j,l}^{\varepsilon}=\{s_{j,l}^{\varepsilon},
s_{i,l}^{-1}s_{k,l}^{-1}s_{i,l}s_{k,l}\} s_{i,k}.\label{E12}
\end{align}
Also for $j<i<k<l$ or $i<k<j<l$, and $\varepsilon, \delta =\pm1$
\begin{align}
&s_{i,k}^{\delta}s_{j,l}^{\varepsilon}=s_{j,l}^{\varepsilon}s_{i,k}^{\delta}.\label{E13}
\end{align}
\end{lemma}

\begin{proof}
 We provide a proof of (\ref{E7}) for
$\varepsilon=1$ as a typical example. We use arguments given by  Markov in ~\cite{Mar45}.
First we can easily see that the relation holds for $j=k-1$,
$l=k+1$. We assume that (\ref{E7}) holds for $j<k<l$. Then we prove
it for $j-1$ and $l+1$ using Lemma \ref{L1}. We deduce the following
equalities by direct computation:
\begin{align*}
&s_{j,k}^{-1}s_{k,l+1}=s_{j,k}^{-1}\{s_{k,l},\sigma_l^{-1}\}
=\{\{s_{k,l},s_{j,l}^{-1}\},\sigma_l^{-1}
\}s_{j,k}^{-1}=\{s_{k,l+1},s_{j,l+1}^{-1}\}s_{j,k}^{-1},\\
&s_{j-1,k}^{-1}s_{k,l}=\{s_{j,k}^{-1},\sigma_{j-1}\}s_{k,l}
=\sigma_{j-1}^{-1}s_{j,k}^{-1}\sigma_{j-1}
s_{k,l}=\sigma_{j-1}^{-1}s_{j,k}^{-1}s_{k,l}\sigma_{j-1}=\\
&\{\{s_{k,l},s_{j,l}^{-1}\}s_{j,k}^{-1},\sigma_{j-1}\}
=\{s_{k,l},s_{j-1,l}^{-1}\}s_{j-1,k}^{-1}.
\end{align*}
This shows that (\ref{E7}) holds.
\end{proof}

Finally, we have the following relations in the braid group $B_n$( see
~\cite[Lemma 5]{Mar45}). A proof is fairly simple.

\begin{lemma}\label{L3}
The following relations hold in the braid group $B_n$:
\begin{align}
&\sigma_{j}^{-1}\sigma_k^{-1}= \sigma_k^{-1}\sigma_{j}^{-1},\quad j<k-1,\label{E14}\\
&\sigma_{j,j+1}\sigma_{k,j+1}=\sigma_{k,j+1}\sigma_{j-1,j},\quad k<j,\label{E15}\\
&\sigma_i^{-2}=s_{i,i+1}^{-1}, \label{E16}\
\end{align}
\end{lemma}

We now call the relations (\ref{E1})--(\ref{E16}) together with the
 trivial relations
$$
s_{i,j}^{\pm 1}s_{i,j}^{\mp 1}=1
$$
the Artin-Markov relations $\mathbf{S}$ for the braid group $B_n$ in terms of the
Artin-Burau generators.

Using the above relations $\mathbf{S}$ together with the definition
$\sigma_i=s_{ii+1}\sigma_i^{-1}$, we can deduce \ori{standard}
Artin's relations for $B_n$.

 Namely, in relation (\ref{E15}) we let $k=j-1$. Then  we
have
$$
\sigma_j^{-1}\sigma_{j-1}^{-1}\sigma_j^{-1}=\sigma_{j-1}^{-1}\sigma_j^{-1}\sigma_{j-1}^{-1}.
$$
Also
$$
\sigma_i^{-1}\sigma_i=\sigma_i^{-1}s_{i,i+1}\sigma_i^{-1}=s_{i,i+1}\sigma_i^{-2}
=s_{i,i+1}s_{i,i+1}^{-1}=1,
$$
and the same for $\sigma_i\sigma_i^{-1}=1$.

\begin{cor}\label{L4} The following relations are deduced by the ELW of $\mathbf{S}$
(to be more precis, by the ELW of  (\ref{E7})-- (\ref{E16})):
\begin{align*}
&\sigma_{i,j}\sigma_{k,j}=\sigma_{k,j}\sigma_{i-1,j-1},\quad k<i,\\
&\sigma_{i,j}\sigma_{k,j}=s_{i,k+1}^{-1}\sigma_{k+1,j}\sigma_{i,j-1},\quad
i\leq k.
\end{align*}
\end{cor}

\section{Main Theorem}
Using the Artin-Markov relations given in the Section 3, we
establish the following theorem.
\begin{thm}\label{T1}
The Artin-Markov relations form a minimal \GS\ basis of the braid
group $B_n$ in term of the Artin-Burau generators relative to the
inverse tower order of words.
\end{thm}

\begin{proof}
There are no inclusion compositions of  defining relations. We only
need to check all possible intersection compositions. Let us do some for examples.

Let us check a composition of intersection of two relations $f,g$ of the form (\ref{E8})
relative to the ambiguity
$$
w=(ij)(jk)(kl),\quad i<j<k<l,
$$
where $(ij)=s_{i,j}$. We have
$$
f=(ij)(jk)-\{(jk),(ik)(jk)\}(ij),\ 
g=(jk)(kl)-\{(kl), (jl)(kl)\}(jk).
$$

We need to prove that
\begin{align}
&\{(jk),(ik)(jk)\}(ij)(kl)\equiv (ij)\{(kl),(jl)(kl)\}(jk) \bmod (\mathbf{S},w).
\label{E17}
\end{align}

In fact, by computation we deduce that
 
\begin{align*}
&(ij)\{(kl),(jl)(kl)\}(jk)\equiv \{(kl),\{(jl),(il)(jl)\}(kl)\}(ij)(jk)\equiv 
\end{align*}
\begin{align}
&\{(kl),\{(jl),(il)(jl)\}(kl)\}\{(jk),(ik)(jk)\}(ij).\label{E18}
\end{align}

For the left hand side of (\ref{E17}) we have 

\begin{align*}
&\{(jk),(ik)(jk)\}(ij)(kl)\equiv\{(jk),(ik)(jk)\}(kl)(ij)\equiv\\
&(jk)^{-1}(ik)^{-1}(jk)(ik)(jk)(kl)(ij)\equiv\\  
&(jk)^{-1}(ik)^{-1}(jk)(ik)\{(kl),(jl)(kl)\}(jk)(ij)\equiv\\
&(jk)^{-1}(ik)^{-1}(jk)\{\{(kl),(il)(kl)\} ,\{(jl),(il)^{-1}(kl)^{-1}(il)(kl) \}\\
&\{(kl),(il)(kl)\}\}(ik)(jk)(ij)\equiv(jk)^{-1}(ik)^{-1}(jk)\{(kl),\\
&\{(jl),(il)^{-1}(kl)^{-1}\}(kl)(il)(kl)\}(ik)(jk)(ij)\equiv\\
&(jk)^{-1}(ik)^{-1}(jk)\{(kl),(il)(jl)(kl)\}(ik)(jk)(ij)\equiv\\
&(jk)^{-1}(ik)^{-1}\{\{(kl),(jl)(kl)\},(il)\{(jl),(kl)\}\{(kl),(jl)(kl)\}\}(jk)(ik)(jk)(ij)\equiv\\
&(jk)^{-1}\{\{\{(kl),(il)^{-1}\},X\},\{(il),(kl)^{-1}(il)^{-1}\}X\}(ik)^{-1}(jk)(ik)(jk)(ij)\equiv
\end{align*}
\begin{align*}
&(X=\{(jl),(kl)(il)(kl)^{-1}(il)^{-1}\}\{(kl),(il)^{-1}\}=\{(kl),(il)^{-1}\}\{(jl),(kl)\})
\end{align*}
\begin{align*}
&(jk)^{-1}\{\{\{(kl),(il)^{-1}\},\{(jl),(kl)\}\},(il)(jl)(kl)\}(ik)^{-1}(jk)(ik)(jk)(ij)\equiv\\
&\{\{\{\{(kl),(jl)^{-1}\},(il)^{-1}\},\{\{(jl),(kl)^{-1}(jl)^{-1}\},\{(kl),(jl)^{-1}\}\}\},\\
&(il)\{(jl),(kl)^{-1}(jl)^{-1}\}\{(kl),(jl)^{-1}\}\}(jk)^{-1}(ik)^{-1}(jk)(ik)(jk)(ij)=
\end{align*}
\begin{align}
&\{\{\{(kl),(jl)^{-1}(il)^{-1}\},(jl)\},(il)(jl)(kl)\}\{(jk),(ik)(jk)\}(ij).
\label{E19}
\end{align}

Words (\ref{E18}) and (\ref{E19}) are the same for

\begin{align*}
&\{(kl),\{(jl),(il)(jl)\}(kl)\}=\{\{\{(kl),(jl)^{-1}(il)^{-1}\},(jl)\},(il)(jl)(kl)\}.
\end{align*}

Thus, (\ref{E17}) is verified and the composition is checked.

To check  the compositions  of relations $f, g$ of the forms
($\ref{E12}), (\ref{E2})$ respectively, we let
$w=\sigma_q^{-1}(jk)(kl)$ and let
\begin{align*}
f&=\sigma_q^{-1}s_{jk}-s_{jk}\sigma_q^{-1},\quad q\neq
j-1,j,k-1,k,\\
g&=(jk)(kl)-\{(kl),(jl)(kl)\}(jk).
\end{align*}
Then we have
$$
(f,g)_w=\sigma_q^{-1}\{(kl),(jl)(kl)\}(jk)-(jk)\sigma_q^{-1}(kl).
$$
If $q\neq l-1,l$ then it is clear that the composition is trivial.
So, we need to consider the following cases:

a) $q=l$. We have
\begin{align*}
\sigma_l^{-1}\{(kl),(j&l)(kl)\}(jk)\equiv\{\{(k,l+1),(l,l+1)\},\{(j,l+1),(l,l+1)\}\\[6pt]
\{(k,l+1),(&l,l+1)\}\}(jk)\sigma_l^{-1}\equiv
\{(k,l+1),(j,l+1)(k,l+1)(l,l+1)\}(jk)\sigma_l^{-1},\\[6pt]
(jk)\sigma_l^{-1}(kl)&\equiv (jk)\{(k,l+1),(l,l+1)\}\sigma_l^{-1}\\
&\equiv\{(k,l+1),(j,l+1)(k,l+1)(l,l+1)\}(jk)\sigma_l^{-1}.
\end{align*}
Hence, the case is verified.

b) $q=l-1$. We have
\begin{align*}
&\sigma_{l-1}^{-1}\{(kl),(jl)(kl)\}(jk)\equiv\{\{(k,l-1),(j,l-1)(k,l-1)\}(jk)\sigma_{l-1}^{-1},\\
&(jk)\sigma_{l-1}^{-1}(kl)\equiv (jk)(kl)\sigma_{l-1}^{-1}\equiv
\{(k,l-1),(j,l-1)(k,l-1)\}(jk)\sigma_{l-1}^{-1}.
\end{align*}
Hence, the case is also verified.

Finally, we need to check the composition of relations (\ref{E16}).
We first let
\begin{align*}
f&=\sigma_j^{-1}\sigma_k^{-1}\dots
\sigma_j^{-1}-\sigma_k^{-1}\dots
\sigma_j^{-1}\sigma_{j-1}^{-1},\quad k<j,\\
g&=\sigma_j^{-1}\sigma_l^{-1}\dots
\sigma_j^{-1}-\sigma_l^{-1}\dots
\sigma_j^{-1}\sigma_{j-1}^{-1},\quad l<j,\\
\text{ and let } w&=\sigma_j^{-1}\sigma_k^{-1}\dots
\sigma_j^{-1}\sigma_l^{-1}\dots \sigma_j^{-1}.
\end{align*}
Then, we have
$$
(f,g)_w=-\sigma_k^{-1}\dots
\sigma_j^{-1}\sigma_{j-1}^{-1}\sigma_l^{-1}\dots \sigma_j^{-1}+
\sigma_j^{-1}\sigma_k^{-1}\dots
\sigma_{j-1}^{-1}\sigma_l^{-1}\dots
\sigma_j^{-1}\sigma_{j-1}^{-1},
$$
We consider the following cases:

a) $l=j-1$.  In this case, by Corollary \ref{L4}, we have
\begin{align*}
&(f,g)_w\equiv -\sigma_k^{-1}\dots \sigma_j^{-1}s_{j-1,j}^{-1} \sigma_j^{-1}+
\sigma_j^{-1}\sigma_k^{-1}\dots
\sigma_{j-2}^{-1}s_{j-1j}^{-1}\sigma_j^{-1}\sigma_{j-1}^{-1};\\[6pt]
&\sigma_k^{-1}\dots \sigma_j^{-1}s_{j-1,j}^{-1} \sigma_j^{-1}\equiv
s_{j,j+1}^{-1}\sigma_k^{-1}\dots \sigma_{j-1}^{-1}s_{jj+1}^{-1}\equiv
s_{j,j+1}^{-1}s_{kj+1}^{-1}\sigma_k^{-1}\dots
\sigma_{j-1}^{-1};\\[6pt]
&\sigma_j^{-1}\sigma_k^{-1}\dots \sigma_{j-2}^{-1}s_{j-1j}^{-1}\sigma_j^{-1}\sigma_{j-1}^{-1}\equiv
\sigma_j^{-1}s_{kj}^{-1}\sigma_k^{-1}\dots \sigma_{j-2}^{-1}\sigma_j^{-1}\sigma_{j-1}^{-1}\\
&\qquad
\equiv\{s_{kj+1}^{-1},s_{jj+1}\}s_{jj+1}^{-1}\sigma_k^{-1}\dots
\sigma_{j-1}^{-1}\equiv
s_{j,j+1}^{-1}s_{kj+1}^{-1}\sigma_k^{-1}\dots \sigma_{j-1}^{-1},
\end{align*}
and the case is done.

b) $l<j-1$. We use again Corollary \ref{L4} to obtain
\begin{align*}
&\sigma_k^{-1}\dots \sigma_j^{-1}\sigma_{j-1}^{-1}\sigma_l^{-1}\dots
\sigma_j^{-1}\equiv \sigma_k^{-1}\dots
\sigma_j^{-1}\sigma_l^{-1}\dots \sigma_j^{-1}\sigma_{j-2}^{-1}\equiv
\sigma_{kj+1}\sigma_{lj+1}\sigma_{j-2}^{-1}\\
&\equiv\sigma_{lj+1}\sigma_{k-1j}\sigma_{j-2}^{-1}\hbox{\;\ (when
$l<k$), or }\equiv s_{kl+1}^{-1}\sigma_{l+1j+1}
\sigma_{kj}\sigma_{j-2}^{-1} \hbox{\;\ (when $k\leq l$)}.
\end{align*}
If $l<k$, then
\begin{multline*}
\sigma_j^{-1}\sigma_k^{-1}\dots
\sigma_{j-1}^{-1}\sigma_l^{-1}\dots
\sigma_j^{-1}\sigma_{j-1}^{-1}\equiv
\sigma_j^{-1}\sigma_{kj}\sigma_{lj}\sigma_j^{-1}
\sigma_{j-1}^{-1}\\
\equiv\sigma_j^{-1}\sigma_{lj}\sigma_{k-1j-1}\sigma_j^{-1}\sigma_{j-1}^{-1}
\equiv \sigma_j^{-1}\sigma_{lj+1}\sigma_{k-1j} \ \equiv
\sigma_{lj+1}\sigma_{k-1j}\sigma_{j-2}^{-1},
\end{multline*}
and we are done. If $l\geq k$,   then $k\leq l<j-1$ and so we have
\begin{multline*}
\sigma_j^{-1}\sigma_k^{-1}\dots
\sigma_{j-1}^{-1}\sigma_l^{-1}\dots
\sigma_j^{-1}\sigma_{j-1}^{-1}\equiv
\sigma_j^{-1}\sigma_{kj}\sigma_{lj}\sigma_j^{-1}
\sigma_{j-1}^{-1}\\
\equiv
\sigma_j^{-1}s_{kl+1}^{-1}\sigma_{l+1j}\sigma_{kj-1}\sigma_j^{-1}\sigma_{j-1}^{-1}
\equiv \sigma_j^{-1}s_{kl+1}^{-1}\sigma_{l+1j+1}\sigma_{kj},
\end{multline*}
and
$$
\sigma_j^{-1}s_{kl+1}^{-1}\sigma_{l+1j+1}\sigma_{kj}\equiv
s_{kl+1}^{-1}\sigma_{l+1j+1} \sigma_{kj}\sigma_{j-2}^{-1},
$$
as desired.
\end{proof}

 Applying  \CD\ Lemma we
obtain:
\begin{cor}The set of $\mathbf{S} $-irreducible words of $B_n$ corresponding to 
the above \GS\ basis $\mathbf{S} $
consists of the words
\begin{align}\label{E20}
&f_nf_{n-1}\dots f_2\sigma_{i_nn}\sigma_{i_{n-1}n-1}\dots \sigma_{i_22},
\end{align}
where $f_j$ are free irreducible words in $\{s_{ij}\mid i< j\},\ 2\leq
j\leq n$.
\end{cor}

\begin{cor}[{Markov-Ivanovsky ~\cite[Theorem 6]{Mar45} and Artin ~\cite[Theorem 17 and
remark of Theorem 18]{Ar47}}] \label{C1} Every word of $B_n$ has a
unique presentation in the form $(\ref{E20})$.
\end{cor}

Let $\Sigma_n$ be the symmetric group, i.e.,
$$
\Sigma_n=\langle s_1,\dots ,s_{n-1} \mid s_i^2=1,
s_{i+1}s_is_{i+1}=s_is_{i+1}s_i, s_is_j=s_js_i,\; i-j>1\rangle,
$$
and let
$$
S_{i,i}=1\hbox{ and }S_{i,j+1}=s_is_{i+1}\dots s_j,\quad i<j.
$$

The following lemma was proved in ~\cite[Theorem 4, Corollary
6]{Mar45}. It also follows from the fact that
$$
\{s_i^2=1,\, s_is_j=s_js_i,\, i-j>1,\,
S_{j,j+1}S_{k,j+1}=S_{k,j+1}S_{j-1,j},\, k<j\}
$$
is a \GS\ basis of $\Sigma_n$ under the deg-inlex order of words
in $\{s_i\}$ (see \cite{BoS01}).

\begin{lemma}\label{L5}
Every element of $\Sigma_n$ has a unique presentation in a form
$$
S_{i_n,n}S_{i_{n-1},n-1}\dots S_{i_2,2},
$$
where $i_j\leq j$ and $2\leq j\leq n$.
\end{lemma}

Let $P_n$ be the group of pure braids. This is the kernel of the natural
homomorphism of $B_n$ onto $\Sigma_n$. From Theorem \ref{T1}, Corollary \ref{C1}
and Lemma \ref{L5} it follows

\begin{cor}[{Markov-Ivanovsky ~\cite[Theorem 8]{Mar45} and Artin ~\cite[Theorem
18]{Ar47}}] $P_n$ is a group with generators $\{s_{ij}\}$ and
defining relations $(\ref{E7})$--$(\ref{E13})$ $($which, together
with the trivial relations, form a minimal \GS\ basis of $P_n$
relative the inverse tower order of words in the generators$)$.
\end{cor}

\end{document}